\documentclass[11pt,reqno,twoside]{amsart}
\usepackage{graphicx}
\usepackage{lmodern} 
\usepackage[T1]{fontenc}
\usepackage[active]{srcltx}
\usepackage{amsmath}
\usepackage{amsthm}
\usepackage{amssymb}
\usepackage{fullpage}
\usepackage[all,cmtip]{xy}
\usepackage{amsfonts}
\usepackage{verbatim}
\usepackage{amscd,latexsym}
\usepackage{tikz}
\usepackage{caption}
\usepackage{array}
\usepackage{multirow}
\usepackage{mathtools}
\usepackage{mathrsfs}
\usepackage{hyperref}
\usepackage{tabularx}
\usepackage{pgfplots}
\pgfplotsset{compat=1.18}
\usepackage{float}
\usepackage{bm}
\usepackage{graphicx}
\usepackage{makecell}
\usepackage{longtable}
\usepackage{mathabx}
\usepackage{mathdots}
\usepackage{mathtools}
\makeatletter
\newcommand{\single}{\let\CS=\@currsize\renewcommand{\baselinestretch}{1.1}\tiny\CS}
\newcommand{\singb}{\let\CS=\@currsize\renewcommand{\baselinestretch}{1}\tiny\CS}
\newcommand{\singa}{\let\CS=\@currsize\renewcommand{\baselinestretch}{1.2}\tiny\CS}
\newcommand{\oneandahalfspacing}{\let\CS=\@currsize\renewcommand{\baselinestretch}{1.5}\tiny\CS}
\newcommand{\singlespacing}{\let\CS=\@currsize\renewcommand{\baselinestretch}{1.6}\large\CS}
\newcommand{\bc}{\begin{center}}
	\newcommand{\ec}{\end{center}}
\newcommand{\be}{\begin{eqnarray}}
	\newcommand{\ee}{\end{eqnarray}}

\newcommand{\Hom}{\operatorname{Hom}}

\newcommand{\Dim}{\operatorname{dim}}

\newcommand{\Irr}{\operatorname{Irr}}
\newcommand{\Ind}{\operatorname{Ind}}

\newcommand{\Rep}{\operatorname{Rep}}

\newcommand{\Span}{\operatorname{Span}}

\newcommand{\beano}{\begin{eqnarray*}}
	\newcommand{\eeano}{\end{eqnarray*}}

\newcommand{\ba}{\begin{array}}
	\newcommand{\ea}{\end{array}}

\makeatletter

\newtheorem{theorem}{Theorem}[section]

\newtheorem{lemma}[theorem]{Lemma}

\theoremstyle{definition}
\newtheorem{definition}{Definition}[section]
\newtheorem{rem}{Remark}[section]
\newtheorem{example}{Example}[section]

\numberwithin{equation}{section}

\DeclareMathOperator{\GL}{GL}
\DeclareMathOperator{\Sp}{Sp}

\DeclareMathOperator{\Par}{P}
\DeclareMathOperator{\G}{G}
\DeclareMathOperator{\F}{F}
\DeclareMathOperator{\N}{N}
\DeclareMathOperator{\Ha}{H}
\DeclareMathOperator{\M}{M}
\DeclareMathOperator{\V}{V}

\allowdisplaybreaks
\begin{document}
	
	\title[Jacquet functor of Symplectic groups ]{On the Jacquet functor of Symplectic groups}
	
	%%%%%%    Information for first author
	
	\author[PREM DAGAR]{PREM DAGAR}
	\address{Department of Mathematics, Indian Institute of Technology Roorkee, Uttarakhand, 247667, India}
	\email{dagarprem5@gmail.com}

	%%%%%%   Information for second author
	
	\author[MAHENDRA KUMAR VERMA]{MAHENDRA KUMAR VERMA}
	\address{Department of Mathematics, Indian Institute of Technology Roorkee, Uttarakhand, 247667, India}
	\email{mahendraverma@ma.iitr.ac.in}


\begin{abstract}
      We prove that, for an equivalence class of irreducible smooth representations of the symplectic group $\Sp_{2n}(\F)$ over a non-Archimedean local field $\F$, the Jacquet functor with respect to the maximal Levi subgroup $\GL_l(\F)\times \Sp_{2(n-l)}(\F)$ is multiplicity-free. The proof is based on an explicit computation of Jacquet modules for a broader family of $\Sp_{2n}(\F)$-representations induced from segments, yielding a detailed structural description that may be of independent interest.

\end{abstract}

\maketitle

\section{Introduction}\label{sec1}
        
         \noindent Let $\G$ be a connected reductive group over a non-Archimedean local field $\F$. A primary focus of representation theory is the study of irreducible representations of groups, as they form the building blocks for understanding the structure of more complex representations. A valuable tool for studying the representations of $\G$ is the Jacquet module. Let us recall the definition of the Jacquet module. Let \(\Par=\M\N\) be a parabolic subgroup of \(\G\), where \(\M\) is a Levi subgroup and \(\N\) its unipotent radical. Denote by \(\mathcal{M}(\G)\) and \(\mathcal{M}(\M)\) the categories of smooth and admissible complex representations of finite length of \(\G\) and \(\M\), respectively. For \((\pi,\V)\in\mathcal{M}(\G)\), define $$\V_{\N}=\V/\Span\{\pi(n)v-v:\ v\in \V,\ n\in \N\}.$$ Since $\M$ normalizes $\N,$ it acts naturally on $\V_{\N}$ by $(r_{\M}^{\G}(\pi))(m)\bigl[v\bigr]=\delta_{\Par}^{-1/2}(m)\,[\pi(m)v] m\in \M,$ where \(\delta_{\Par}\) denotes the modular function of \(\Par\), and \([v]\) is the image of \(v\) in \(\V_{\N}\). The resulting representation \((r_{\M}^{\G}(\pi),\V_{\N})\in\mathcal{M}(\M)\) is called the normalized Jacquet module of \((\pi,\V)\). The functor $$r_{\M}^{\G}:\mathcal{M}(\G)\to\mathcal{M}(\M),\qquad \pi\mapsto r_{\M}^{\G}(\pi),$$ is called the Jacquet functor, or Jacquet restriction functor. The Jacquet functor is exact and left adjoint to parabolic induction i.e., Frobenius reciprocity holds. More precisely, for every smooth representation \(\pi\) of \(\G\) and every smooth representation \(\sigma\) of \(\M\), there is a natural isomorphism  $\Hom_{\G}\!\bigl(\pi,\Ind_{\Par}^{\G}(\sigma)\bigr)\cong\Hom_{\M}\!\bigl(r_{\M}^{\G}(\pi),\sigma\bigr),$ where \(\Ind_{\Par}^{\G}(\sigma)\) denotes the representation of \(\G\) parabolically induced from \(\sigma\). We consider normalized Jacquet and induction functors, i.e., they carry unitarizable representations to unitarizable ones. The Frobenius reciprocity indicates the importance of the Jacquet module of a representation.

         Tadi$\acute{c}$ \cite{MR1658535} proved that the Jacquet module is a useful tool for establishing the irreducibility of an induced representation. To explore this, let $\Par_0 = \M_0\N_0$ be a parabolic subgroup of $\G$ and let $\sigma_0$ be an irreducible admissible representation of $\M_0$. Let $\Par'=\M'\N', \Par''=\M''\N'',\Par'''=\M'''\N'''$ be parabolic subgroups of $\G$ such that $\Par'\subset \Par'', \Par' \subset \Par'''$ and $r_{\M'}^{\G}(\Ind_{\Par_0}^{\G}(\sigma_0)) \neq 0$. Suppose that there exists an irreducible subquotient $\tau''$ of $r_{\M''}^{\G}(\Ind_{\Par_0}^{\G}(\sigma_0)) $ such that for any irreducible subquotient $\tau'''$ of $r_{\M'''}^{\G}(\Ind_{\Par_0}^{\G}(\sigma_0))$, we have $$r_{\M'}^{\M''}(\tau'')+ r_{\M'}^{\M'''}(\tau''')\nleq r_{\M'}^{\G}(\Ind_{\Par_0}^{\G}(\sigma_0)).$$ Then, $\Ind_{\Par_0}^{\G}( \sigma_0)$ is irreducible. Moreover, for reducible induced representations, the multiplicity-freeness can often be determined by analyzing the multiplicity-freeness of their corresponding Jacquet modules. Thus, multiplicity-free Jacquet modules play a central role in this analysis. This leads naturally to the question of which irreducible representations have multiplicity-free Jacquet modules. In this article, we provide a class of irreducible representations of the symplectic group whose Jacquet module, with respect to the maximal Levi subgroup is multiplicity-free.
        
        % Let $\Par_0 = \M_0\N_0$ be a parabolic subgroup of $\G$ and $\sigma_0$ be an irreducible unitarizable admissible representation of $\M_0$. Tadi$\acute{c}$\cite{MR1658535} explored that if the multiplicity of $\sigma_0$ in $r_{\M_0}^{\G}(\Ind_{\Par_0}^{\G}( \sigma_0))$ is one, then $\Ind_{\Par_0}^{\G}( \sigma_0)$ is irreducible. If the multiplicity of  $\sigma_0$ in $r_{\M_0}^{\G}(\Ind_{\Par_0}^{\G}( \sigma_0))$ is two, then $\Ind_{\Par_0}^{\G}( \sigma_0)$ is either irreducible or a direct sum of two irreducible non-isomorphic representations. We cannot assert anything regarding the converse of the previous statements. An irreducible representation may have a Jacquet module with an irreducible subquotient that appears multiple times.
\subsection{Main result}

       Let $\GL_n=\GL_n(\F)$ denote the general linear group over a non-Archimedean local field $\F$, and let $\Sp_n=\Sp_{2n}(\F)$ denote the symplectic group (see \S\ref{symp} for the definition). Let $\Par$ be a maximal parabolic subgroup of $\Sp_n$ with Levi decomposition $\Par=\M\N$, where $\M\simeq\GL_l\times \Sp_{n-l}$ is the corresponding maximal Levi subgroup and $\N$ is its unipotent radical. For the group $\GL_n$, Zelevinsky \cite{z80} computed the Jacquet modules of representations of the form $Z(\Delta)$ (see \S\ref{preln} for the relevant notion). Later, Kret and Lapid \cite{kret} studied the Jacquet modules of ladder representations of $\GL_n$. In \cite{prem}, we studied the Jacquet modules of induced representations of the form $Z(\Delta_1)\times \cdots \times Z(\Delta_r)$ for $\GL_n$. For $\Sp_n$, Mati\'c and Tadi\'c~\cite{tadmat} studied the Jacquet modules of representations of the form $\delta(\Delta)\rtimes \sigma$, where $\delta(\Delta)$ is a square-integrable representation of a general linear group of smaller rank and $\sigma$ is an irreducible supercuspidal representation of a symplectic group of smaller rank. Subsequently, Bo\v{s}njak \cite{bar} studied the Jacquet modules of representations induced from ladder and cuspidal representations. For further developments concerning the structure and properties of Jacquet modules, we refer the reader to \cite{tadic1995structure, aizenbud2012multiplicity, MR1758232, MR3163355, MR3384985, MR3783564}.

        In this article, we study the Jacquet module of the representations of $\Sp_{n}$ of the form $Z( \Delta_1)\times \cdots\times Z( \Delta_r)\rtimes\sigma$, where $Z( \Delta_1)\times \cdots\times Z( \Delta_r)$ is a representation of general linear group of smaller rank and $\sigma$ is an irreducible supercuspidal representation of smaller rank symplectic group. We conclude with a precise statement on the multiplicities appearing in these Jacquet modules.

\begin{definition}
        Let $\Ha$ and $\Ha'$ be two reductive groups over a field $\F$, and let $\mathcal{M}(\Ha)$ and $\mathcal{M}(\Ha')$ represents the category of all smooth representations of $\Ha$ and $\Ha'$, respectively.~We call a functor $\mathcal{F}:\mathcal{M}(\Ha)\rightarrow \mathcal{M}(\Ha')$, multiplicity free if for any irreducible representation $\pi\in \mathcal{M}(\Ha)$, the representation $\mathcal{F}(\pi)$ is multiplicity free representation, i.e., $\dim \Hom_{\Ha'}(\mathcal{F}(\pi), \tau )\leq 1,$  for any irreducible representation $\tau\in \mathcal{M}(\Ha').$
\end{definition}

       In this article, we consider the Jacquet functor $r_{\M}^{\Sp_n}:\mathcal{M}(\Sp_n)\to \mathcal{M}(\M),$ where \(\M\) is a maximal Levi subgroup of \(\Sp_n\). We prove that, when restricted to \(\mathcal{M}_{\Irr}(\Sp_n)\) (see \S\ref{class2} for notation), the functor \(r_{\M}^{\Sp_n}\) is multiplicity-free. Our main result is the following theorem.
       
\begin{theorem}\label{main} 

        Let $\M$ be a maximal levi subgroup of~ $\Sp_n$ and $r_{\M}^{\Sp_{n}}:\mathcal{M}_{\Irr}(\Sp_{n})\rightarrow \mathcal{M}(\M)$ be the Jacquet functor. Then $$ \Dim_{\mathbb{C}}\Hom_{\M}(r_{\M}^{\Sp_{n}}(\pi), \rho) \leq 1,$$	for  $\pi\in\mathcal{M}_{\Irr}(\Sp_{n})$ and any irreducible representation $\rho\in\mathcal{M}(\M).$ 
        	
\end{theorem}

         This work is inspired by the work of Aizenbud and Gourevitch \cite[Theorem A]{aizenbud2012multiplicity}, which proves that the Jacquet functor $r_{\GL_n\times\GL_k}^{\GL_{n+k}} :\mathcal{M}(\GL_{n+k}) \rightarrow \mathcal{M}(\GL_n\times\GL_k)$ is multiplicity-free in the cases \(k=1\) and \(k=2\). In our previous work \cite[Theorem 1.1]{prem}, we studied the Jacquet modules of certain representations of \(\GL_n\) relative to an arbitrary maximal Levi subgroup \(\M\). We also established that, for a distinguished subclass of irreducible representations, the corresponding Jacquet modules are multiplicity-free.

\subsection{Remark} 
        This result is more or less an observation that the Jacquet module of any $\pi \in \mathcal{M}_{\Irr}(\Sp_{n})$ is multiplicity free. While broader classes of representations may yield similar results, our focus lies in establishing this property specifically for the class $\mathcal{M}_{\Irr}(\Sp_{n})$.

\subsection{Organization:}
       Our discussion of the results unfolds in a structured manner. Section 2 gives essential preliminaries and the notations needed for our discourse. Following this, Section 3 deals with the computation of Jacquet modules for representations within $\pi \in \mathcal{M}(\Sp_{n})$ followed by proof of Theorem \ref{main}. To enhance understanding, we conclude this section with an illustrative example.

\section{Preliminaries}

\subsection{Generalities}
This section introduces the notation and preliminary concepts used throughout
the paper. The primary reference is \cite{z80}.

\subsubsection{}
       For a reductive group $\G$, let $\Rep(\G)$ denote the category of smooth admissible complex representations of $\G$, $\Irr(\G)$ the set of equivalence classes of irreducible representations, $\mathcal{C}(\G)$ the set of supercuspidal representations, and $\mathcal{C}_{\Irr}(\G)$ the set of irreducible supercuspidal representations. For $\pi\in \Rep(\G)$, we write $\widetilde{\pi}$ for its contragredient. Throughout this article, all representations are assumed to be smooth and admissible.    
\subsubsection{}\label{symp} 
       Consider the skew-symmetric matrix ~$ \begin{pmatrix}
	0 & -J_n\\
	J_n & 0 \end{pmatrix},$ where $J_n=\left(\delta_{i,n+1-j}\right)_{1\le i,j\le n}.$ The symplectic group $\Sp_{2n}(\F)$ is defined as follows:
% $$   \vspace{0.1cm} $$J_n =	\begin{pmatrix}
% 	& & & & 1  \\
% 	& &  & 1&  \\
% 	&  & \iddots  & &  \\
% 	&1 & &  &\\
	
%        \end{pmatrix}.$$\vspace{0.2cm} 
 $$\Sp_{2n}(\F)=\bigg\{g\in \GL_{2n}(\F): g^t\begin{pmatrix}
	0 & -J_n\\
	J_n & 0
       \end{pmatrix} g=\begin{pmatrix}
	0 & -J_n\\
	J_n & 0
        \end{pmatrix} \bigg\}.$$ 
      Throughout this paper, we denote $\GL_n(\F)$ by $\GL_n$ and $\Sp_{2n}(\F)$ by $\Sp_n$. Let $\alpha=(n_1,\ldots,n_r)$ be a partition of $m\leq n$. Given such a partition $\alpha $, the set $\M_\alpha \simeq \GL_{n_1}\times\cdots\times \GL_{n_r}\times \Sp_{n-m}$ is a Levi subgroup of~ $\Sp_{n}$. In particular, if $\alpha=(m),$ then  $\M_\alpha\simeq\GL_m\times\Sp_{n-m}$ is the maximal Levi subgroup of $\Sp_{n}$.
\subsubsection{}\label{preln}

       Let $\nu$ be the character of $\GL_n$ defined by $\nu(g)=|\det(g)|$. Let $\wp$ denote the set of equivalence classes of irreducible cuspidal representations of the groups $\GL_n$, $n = 1,2,\dots$. For $\rho \in \wp$ and integers $a \le b$, the set $\Delta = [a,b]_{(\rho)} = \{\nu^{a}\rho, \nu^{a+1}\rho, \ldots, \nu^{b}\rho\}$ is called a segment (associated to $\rho$). The integers \(a\) and \(b\) are called the beginning and the end of \(\Delta\), respectively, and $l(\Delta)=b-a+1$ is called the length of \(\Delta\). The contragredient segment of \(\Delta\) is defined by $\widetilde{\Delta}=[-b,-a]_{(\widetilde{\rho})}.$
        
\begin{definition}
         Let $\Delta_1 = [a_1,b_1]_{(\rho)} $ and $\Delta_2 = [a_2,b_2]_{(\rho)} $ be two segments. We say that $\Delta_1$ and $\Delta_2$ are linked if $\Delta_1 \nsubseteq \Delta_2$, $\Delta_2 \nsubseteq \Delta_1$, and $\Delta_1 \cup \Delta_2$ is also a segment.  We say that $\Delta_1$ precedes $\Delta_2$ if $a_1<a_2$, $b_1<b_2$, and $a_2 \le b_1+1$. 
\end{definition}

     We have a straightforward lemma regarding the non-linkedness of two smaller segments.
\begin{lemma}\label{remark}
	Let $\Delta_1 = [a_1,b_1]_{(\rho)}$ and $\Delta_2 = [a_2,b_2]_{(\rho)}$ be non-linked segments such that $\Delta_1\not\subset\Delta_2$ and $\Delta_2\not\subset\Delta_1$. Then the segments $ [{a_1}, {c_1}]_{(\rho)}$  and $ [{a_2}, {c_2}]_{(\rho)}$ are also not linked for any $a_i\leq c_i\leq b_i,~i=1,2.$ Also, the segments $ [{c_1}, {b_1}]_{(\rho)}$  and $ [{c_2}, {b_2}]_{(\rho)}$ are not linked.
\end{lemma}

%      The following proposition by Casselman \cite{casselman1995introduction} provides a result on the irreducibility of the representation of the group $\Ha_1\times\Ha_2$.
% \begin{proposition}
%       Let $\pi_i$ be smooth irreducible representations of groups $\Ha_i,i = 1, 2.$ Then $\pi_1\otimes\pi_2$ is a smooth irreducible representation of group $\Ha_1\times\Ha_2$.
% \end{proposition}

      Segments provide a systematic way to construct the representations of reductive groups, forming the basis of the Zelevinsky and Langlands classifications. To each segment $\Delta = [a,b]_{(\rho)}$, we associate an irreducible representation $Z(\Delta)$ (resp. $\delta(\Delta)$), defined as the unique irreducible subrepresentation (resp. quotient) of $\nu_\rho^{a}\rho \times \nu_\rho^{a+1}\rho \times \cdots \times \nu_\rho^{b}\rho$. The Zelevinsky classification \cite{z80} states that every irreducible representation of $\GL_n$ is the unique irreducible subrepresentation of an induced representation of the form $Z(\Delta_1)\times \cdots \times Z(\Delta_r).$ We recall Zelevinsky irreducibility criterion (\cite[Theorem 4.2]{z80}) for such induced representations.

\begin{theorem}\label{linked}
	Let $\Delta_1,\ldots,\Delta_r$ be the segments of irreducible cuspidal
representations. Then following are
	equivalent:
	\begin{enumerate}
		\item[\upshape(1)] The representation $Z( \Delta_1)\times \cdots\times Z( \Delta_r) $ is irreducible.
		\item[\upshape(2)] 	 For each $i, j= 1, \ldots, r,$ the segments $\Delta_i$ and $\Delta_j$ are not linked.
	\end{enumerate}
\end{theorem}

       The following theorem analyzes the isomorphism between two induced representations of type $Z( \Delta_1)\times \cdots \times Z( \Delta_r) $, given by Zelevinsky \cite[Theorem 6.4]{z80}.

\begin{theorem}\label{isom}
	   Let $(\Delta_1,\ldots,\Delta_r)$ and $(\Delta_1',\ldots,\Delta_r')$ be ordered sequences of segments in $\wp$. Suppose that any of the following conditions holds:
	\begin{enumerate}
		\item[\upshape(1)] $(\Delta_1,\ldots,\Delta_r)$ differs from $(\Delta_1',\ldots,\Delta_r') $ only by a transposition of two neighbors which
		are not linked.
		\item[\upshape(2)] Both  $(\Delta_1,\ldots,\Delta_r)$  and $(\Delta_1',\ldots,\Delta_r') $ satisfy the condition: For each pair of indicies $i,j$ such that $i<j$, $\Delta_i$ does not precedes $\Delta_j $ (same for $\Delta_i'$), and are equal up to a rearrangement.
	\end{enumerate}	
	Then $Z( \Delta_1)\times \cdots\times Z( \Delta_r) \simeq Z(\Delta_1')\times \cdots\times Z( \Delta_r') $.
\end{theorem}

\subsubsection{}
        % Let $(n_1,\ldots,n_r)$ be a partition of $m\leq n$. Let $\Delta_i= [{a_i}\rho_i,{b_i}\rho_i]$ be the  segments, where $\rho_i$ is an  irreducible supercuspidal representation of $\GL_{m_i}$,  for $i=1,2,\ldots,r.$  Let $b_i-a_i+1=k_i$ and $k_im_i=n_i$. Therefore, $Z( \Delta_i)$ is an irreducible representation of $\GL_{n_i}$ and $Z( \Delta_1)\times\cdots\times Z( \Delta_r) $ be a representation of $\GL_m$. Let $\sigma$ be an irreducible supercuspidal representation of $\Sp_{n-m}$.
        
        It is well known that every irreducible representation of $\Sp_n$ can be realized as a subquotient of $Z(\Delta_1) \times \cdots \times Z(\Delta_r) \rtimes \sigma,$ where $Z(\Delta_1) \times \cdots \times Z(\Delta_r)$ is a representation of a general linear group of smaller rank, and $\sigma$ is an irreducible supercuspidal representation of a symplectic group of smaller rank. From now on, we fix an
irreducible supercuspidal representation $\sigma$ of $\Sp_r$ for $0\leq r\leq n.$ Mui$\acute{c}$ \cite[Theorem 0.1]{muic} established following theorem regarding the irreducibility of $Z( \Delta_1)\times \cdots\times Z( \Delta_r)\rtimes\sigma$.
        
\begin{theorem}
        The representation $Z( \Delta_1)\times \cdots\times Z( \Delta_r)\rtimes\sigma$ is irreducible if and only if following holds:
        \begin{itemize}
       \item[\upshape(1)] $Z( \Delta_i)\times Z( \Delta_j)$ is irreducible for all $i \neq j$,
        \item[\upshape(2)] $Z( \Delta_i)\times Z( \widetilde{\Delta_j)}$ is irreducible for all $i \neq j$,
        \item[\upshape(3)] $Z( \Delta_i)\rtimes \sigma$ is irreducible for all $i.$
        \end{itemize}
\end{theorem}
       The above theorem was originally proved for the representation $\delta( \Delta_1)\times \cdots\times \delta( \Delta_r)\rtimes\sigma$. However, since $\delta(\Delta)$ and $Z(\Delta)$ are dual of each other, the theorem also applies similarly to the representation $Z( \Delta_1)\times \cdots\times Z( \Delta_r)\rtimes\sigma$. 
       
%        Now, the question is about the irreducibility of the representation $Z( \Delta_i)\rtimes \sigma$. The following theorem by Lapid and Tadi$\acute{c}$ \cite[Theorem 3.10]{MR4084162} gives the irreducibility criteria of $Z( \Delta_i)\rtimes \sigma$.
% \begin{proposition}
%           Let $Z( \Delta)$ is an irreducible representation of $\GL_n$ and $\sigma$ be an irreducible cuspidal representation of $\Sp_m$. Then, the following are equivalent.
%           \begin{enumerate}
%               \item $Z( \Delta)\rtimes \sigma$ is irreducible.
%               \item $\delta( \Delta)\rtimes \sigma$ is irreducible.
%               \item $\rho\rtimes \sigma$ is irreducible for every $\rho\in\Delta$.
%           \end{enumerate}
% \end{proposition}
%        An immediate corollary follows from this result, which will be helpful in the main proof of Theorem $\ref{main}$.
% \begin{corollary}\label{2.5}
%         If $\Delta'\subset\Delta$ and $Z( \Delta)\rtimes \sigma$ is irreducible, then $Z( \Delta')\rtimes \sigma$ is also irreducible.
% \end{corollary} 
% \begin{proof}
%         If $Z( \Delta)\rtimes \sigma$ is irreducible, then by above theorem $\rho\rtimes \sigma$ is irreducible for every $\rho\in\Delta$. Consequently, $\rho\rtimes \sigma$ is irreducible for every $\rho\in\Delta'$, which implies that $Z( \Delta')\rtimes~ \sigma$ is irreducible.
% \end{proof}

\subsubsection{}\label{class2}
       It follows from the Theorem \ref{linked} that for the representation  $Z( \Delta_1)\times \cdots\times Z( \Delta_r)\rtimes\sigma$ to be irreducible, it is necessary that $\Delta_i$ and $\Delta_j,\widetilde{\Delta_j}$ be unlinked for all $i\neq j.$

\noindent   
      Let $\mathcal{M}_{\Irr}(\Sp_{n})$ denote the equivalence class of irreducible representations of the form $Z( \Delta_1) \times \cdots \times Z( \Delta_r) \rtimes \sigma$, where $\Delta_i \not\subset \Delta_j$ or $\widetilde{\Delta_j}$ for $i \neq j$, and $\rho_i \not\cong \widetilde{\rho_i}$ for each $i$. Note that \(\mathcal{M}_{\Irr}(\Sp_n)\) is the class of representations for which we aim to prove that the Jacquet functor is multiplicity-free.
% \begin{itemize}

%      \item $\Delta_i$ and $\Delta_j,\widetilde{\Delta_j}$ unlinked for all $i\neq j,$
%     \item $\Delta_i\not\subset\Delta_j, \widetilde{\Delta_j} $ for $i\neq j$.
%     \item $\rho_i\ncong \widetilde{\rho_i}$ for each $i$.
% \end{itemize}
 % The above conditions are imposed to ensure the irreducibility of the representations and to obtain irreducible subquotients in the Jacquet module.

%\subsubsection{}\label{class3}
% For a partition $\alpha=(m),$ let $\mathcal{M}(\M_{\alpha})$ denote the class of representations of the form $\bigoplus_{i,j}\rho_i\otimes\rho_j$, where $\rho_i\in \mathcal{M}(\GL_m)$, $\rho_j\in\mathcal{M}(\Sp_{n-m})$ and $i,j$ ranges over finite natural numbers.

\subsubsection{} 
        Let $\mathcal{R}(\GL_n)$ denote the Grothendieck group of smooth and admissible complex representations of the finite length of $\GL_n$ and  $\mathcal{R}(\GL)=\bigoplus_{n\geq 0}\mathcal{R}(\GL_n).$ Similarly, $\mathcal{R}(\Sp)=\bigoplus_{n\geq 0}\mathcal{R}(\Sp_{n}).$ 
 
% To state the next result, we now introduce the following mappings.
\subsubsection{}
        Given $\sigma \in \mathcal{M}(\Sp_{n})$ and Levi subgroup $\M_\alpha$, we will denote $r_{\M_\alpha}^{\Sp_{n}}(\sigma)$  by $r_{\alpha}(\sigma)$, whenever required. In the absence of any ambiguity in $\alpha$, we may express $\mathcal{M}(\M_{\alpha})$  simply by $\mathcal{M}(\M).$ Consider the maps  $$\mu^*: \mathcal{R}(\Sp_{n}) \rightarrow \mathcal{R}(\GL) \otimes \mathcal{R}(\Sp)$$ 
     \begin{equation}\label{map1}
	\sigma\mapsto 1\otimes\sigma+ \sum_{m=1}^{n}r_{(m)}(\sigma).
     \end{equation}
        and  $$m^*: \mathcal{R}(\GL_n) 
        \rightarrow \mathcal{R}(\GL) \otimes \mathcal{R}(\GL)$$
\begin{equation}\label{map2}
	\pi \mapsto \sum_{m=0}^{n}r_{(m,n-m),(n)}(\pi).
\end{equation}

        \noindent where  $r_{(m,n-m),(n)}(\pi)$ is the normalized Jacquet module of $\pi$ with respect to maximal Levi $\GL_m\times\GL_{n-m}.$ The mappings $\mu^*$ and $m^*$ extend to  $\mathcal{R}(\Sp)$ and $\mathcal{R}(\GL)$, respectively. Consequently, $\mathcal{R}(\Sp)$ and $\mathcal{R}(\GL)$ forms a Hopf algebra structure. 
\subsubsection{}
        Similarly define the induction map $$m: \mathcal{R}(\GL) 
        \otimes \mathcal{R}(\GL) \rightarrow \mathcal{R}(\GL)$$ 
   \begin{equation}\label{mmap}
	\pi_1\otimes\pi_2\mapsto\pi_1\times \pi_2,
    \end{equation} and $$s: \mathcal{R}(\GL) \otimes \mathcal{R}(\GL) \rightarrow \mathcal{R}(\GL)\otimes\mathcal{R}(\GL)$$ 
    \begin{equation}\label{smap}
	\sum_i x_i\otimes y_i \mapsto \sum_i y_i\otimes x_i.
   \end{equation}
\subsubsection{}
        \noindent Now, we define the composition map $$M^*: \mathcal{R}(\GL) \rightarrow \mathcal{R}(\GL) \otimes \mathcal{R}(\GL)$$ by \begin{equation}\label{map3}
	M^*=(m\otimes1)\circ(\widetilde{.} \otimes m^*)\circ s\circ m^*,
\end{equation}

           Now, we state a theorem, given by Tadi$\acute{c}$ \cite[Theorem 5.2]{tadic1995structure}, which we use to compute the Jacquet module of all representations of the form $Z( \Delta_1) \times \cdots \times Z( \Delta_r)\rtimes\sigma.$
\begin{theorem}\label{main1}
	    For $\pi\in \mathcal{R}(\GL)$ and $\sigma\in \mathcal{R}(\Sp)$ $$\mu^*(\pi\rtimes\sigma)=M^*(\pi)\rtimes\mu^*(\sigma),$$ where~the right hand side is identified by
	$$(\pi_1\otimes\pi_2)\rtimes(\pi'\otimes\sigma')=(\pi_1\times\pi')\otimes(\pi_2\rtimes\sigma').$$
\end{theorem}
        Our primary objective is to compute the Jacquet module of all representations of type $Z( \Delta_1) \times \cdots \times Z( \Delta_r) \rtimes \sigma$. For this, we apply Theorem \ref{main1}, which reduces the problem to computing $M^*(Z( \Delta_1) \times \cdots \times Z( \Delta_r))$ and $\mu^*(\sigma)$. Since $\sigma$ is supercuspidal,  $\mu^*(\sigma)=1\otimes\sigma$. Furthermore, by Equation (\ref{map3}), computing $M^*(Z( \Delta_1) \times \cdots \times Z( \Delta_r))$ requires us to first calculate $m^*(Z( \Delta_1) \times \cdots \times Z( \Delta_r))$, which we recall in the next subsection.

\subsection{Jacquet module of representation of general linear group} 
        Now, we recall the calculation of the Jacquet module of representations $Z( \Delta_1)\times\cdots\times Z( \Delta_r)$ of $\GL_m$, detailed in our previous work \cite[ Section 3.2]{prem}. Let $\beta=(n_1,\ldots,n_r)$ and $\gamma=(q,m-q)$ be partitions of $m$. For each $i=1,\ldots,r$, let $\Delta_i=[a_i,b_i]_{(\rho_i)}$ be a segment with $\rho_i\in\mathcal{C}_{\mathrm{Irr}}(\GL_{m_i})$, and set $k_i=b_i-a_i+1$ so that $n_i=k_i m_i$. Then $Z(\Delta_i)\in\Irr(\GL_{n_i})$, and $Z(\Delta_1)\times\cdots\times Z(\Delta_r)$ is a representation of $\GL_m$. Let $M^{\beta,\gamma}$ be the set of matrices $W=(b_{ij})$, where $(b_{ij})$ is the matrix of order $r\times2$ with entries $(b_{i1})=l_i$ and $(b_{i2})=q_i$ satisfying $l_i+q_i=n_i,~ \sum_{i=1}^r l_i=q,~\sum_{i=1}^r q_i=m-q.$ Then, we have
          \begin{align}\label{jacGL}
r_{(q,m-q),(m)}(Z(\Delta_1) \times \cdots \times Z(\Delta_r))
&= \sum_{W \in M^{\beta,\gamma}} \Biggl(
\Bigl( Z([a_1, p_1^{(w)}-1]_{(\rho_1)}) \times \cdots \times Z([a_r, p_r^{(w)}-1]_{(\rho_r)}) \Bigr) \notag \\
&\otimes \Bigl( Z([p_1^{(w)}, b_1]_{(\rho_1)}) \times \cdots \times Z([p_r^{(w)}, b_r]_{(\rho_r)}) \Bigr)
\Biggr) \notag \\
&= \sum_{W \in M^{\beta,\gamma}} \sigma_1^{(w)} \otimes \sigma_2^{(w)} \quad (\text{say})
\end{align}

    such that the non-negative integer $p_i^{(w)}$ satisfies:

\begin{itemize}
	\item $ a_i-1\leq p_i^{(w)}\leq b_i+1$
	\item $l_i=m_i(p_i^{(w)}-a_i)$
	\item $q_i=m_i(b_i-p_i^{(w)}+1)$.
\end{itemize}

\begin{rem}\label{rem}
        Observe that if we impose the condition that $\Delta_i$ and $\Delta_j$ are non-linked with $\Delta_i\not\subset\Delta_j$ for $i\neq j$ and $i,j=1,2,\ldots,r$, then $\sigma_1^{(w)}$ and $\sigma_2^{(w)}$ are irreducible representations of $\GL_q$ and $\GL_{m-q}$, respectively. Also the representations $\sigma_i^{(w)}$ for $i=1,2;$ depends on $W\in M^{\beta,\gamma}$ and all matrices $W$  are distinct. This implies that all the factors in  $r_{(q,m-q),(m)}(Z( \Delta_1)\times\cdots\times Z( \Delta_r))$ are distinct and has no repetitions. 
\end{rem}

\section{Jacquet module of representation of Symplectic group}

In this section, we compute the Jacquet module of $\pi=Z( \Delta_1)\times\cdots\times Z( \Delta_r)\rtimes\sigma \in \mathcal{M}(\Sp_{n}).$ Let $(n_1,\ldots,n_r)$ be a partition of $m\leq n$. Let $\Delta_i= [a_i,b_i]_{(\rho_i)}$ be the  segments, where $\rho_i\in \mathcal{C}_{\Irr}(\GL_{m_i})$, for $i=1,2,\ldots,r.$  Let $b_i-a_i+1=k_i$ and $k_im_i=n_i$. Therefore, $Z( \Delta_i)$ is an irreducible representation of $\GL_{n_i}$ and $Z( \Delta_1)\times\cdots\times Z( \Delta_r) $ is a representation of $\GL_m$. Let $\sigma\in\mathcal{C}(\Sp_{n-m})$, and $\M=\GL_l\times\Sp_{n-l}$ be a maximal Levi subgroup of $\Sp_{n}$.~Applying  Theorem $\ref{main1}$, we first compute $M^*(Z( \Delta_1)\times\cdots\times Z( \Delta_r))$. By (\ref{map3}) and (\ref{jacGL}) (substituting $p_i^{(w)}$ with $p_i^{(j,q)}$ for notational convenience), we obtain

\bigskip

\noindent\underline{\textbf{Step 1:} Applying $m^*$}

\begin{align}
Z(\Delta_1) \times \cdots \times Z(\Delta_r)
&\xmapsto{\quad m^* \quad}
\sum_{q=0}^{m} \sum_{j=1}^{j_q}
\underbrace{\Bigl( Z\bigl([a_1, p_1^{(j,q)}-1]_{(\rho_1)}\bigr) \times \cdots \times
Z\bigl([a_r, p_r^{(j,q)}-1]_{(\rho_r)}\bigr) \Bigr)}_{A_{j,q}} \notag \\
&\qquad \hspace{2cm} \otimes
\underbrace{\Bigl( Z\bigl([p_1^{(j,q)}, b_1]_{(\rho_1)}\bigr) \times \cdots \times
Z\bigl([p_r^{(j,q)}, b_r]_{(\rho_r)}\bigr) \Bigr)}_{B_{j,q}}.\notag
\end{align}
$A_{j,q}$ and $B_{j,q}$ are representations of $\GL_q$ and $\GL_{m-q}$, respectively. Here $q$ varies from $0$ to $m$ and non-negative integer $p_i^{(j,q)}$ satisfies:
$$
\sum_{i=1}^{r} m_i\bigl(p_i^{(j,q)} - a_i\bigr) = q
\quad \text{and} \quad
\sum_{i=1}^{r} m_i\bigl(b_i - p_i^{(j,q)} + 1\bigr) = m - q.
$$

\noindent\underline{\textbf{Step 2:} Applying $s$}
\begin{align}
\sum_{q=0}^{m} \sum_{j=1}^{j_q} A_{j,q} \otimes B_{j,q}
&\xmapsto{\quad s \quad}
\sum_{q=0}^{m} \sum_{j=1}^{j_q} B_{j,q} \otimes A_{j,q}. \notag 
\end{align}

\noindent\underline{\textbf{Step 3:} Applying $\widetilde{.} \otimes m^*$}

\noindent Following~\eqref{map3}, we apply $\widetilde{.}$ to $B_{j,q}$ and $m^*$ to $A_{j,q}$. Write
$$
m^*(A_{j,q}) = \sum_{s=0}^{q} \sum_{u=1}^{u_s(j,q)} A'_{u,s,j,q} \otimes A''_{u,s,j,q},
$$
where
\begin{align*}
A'_{u,s,j,q} &= Z\bigl([a_1, p_1^{(u_s(j,q))}-1]_{(\rho_1)}\bigr) \times \cdots \times
Z\bigl([a_r, p_r^{(u_s(j,q))}-1]_{(\rho_r)}\bigr), \\
A''_{u,s,j,q} &= Z\bigl([p_1^{(u_s(j,q))}, p_1^{(j,q)}-1]_{(\rho_1)}\bigr) \times \cdots \times
Z\bigl([p_r^{(u_s(j,q))}, p_r^{(j,q)}-1]_{(\rho_r)}\bigr),
\end{align*}
and denote
\[
\widetilde{B}_{j,q} =
Z\bigl([-b_1, -p_1^{(j,q)}]_{(\tilde{\rho}_1)}\bigr) \times \cdots \times
Z\bigl([-b_r, -p_r^{(j,q)}]_{(\tilde{\rho}_r)}\bigr).
\]
Then
\begin{align}
\sum_{q=0}^{m} \sum_{j=1}^{j_q} B_{j,q} \otimes A_{j,q}
&\xmapsto{\;\widetilde{.}\,\otimes\, m^*\;}
\sum_{q=0}^{m} \sum_{j=1}^{j_q} \sum_{s=0}^{q} \sum_{u=1}^{u_s(j,q)}
\widetilde{B}_{j,q} \otimes A'_{u,s,j,q} \otimes A''_{u,s,j,q}. \notag
\end{align}
The resulting components $A'_{u,s,j,q}$ and $A''_{u,s,j,q}$ are representations of $\GL_s$ and $\GL_{q-s}$, respectively. Here $s$ varies from $0$ to $q$ and non-negative integer $p_i^{(u_s(j,q))}$ satisfies:
$$
\sum_{i=1}^{r} m_i\bigl(p_i^{(u_s(j,q))} - a_i\bigr) = s
\quad \text{and} \quad
\sum_{i=1}^{r} m_i\bigl(p_i^{(j,q)} - p_i^{(u_s(j,q))}\bigr) = q - s.
$$

\bigskip

\noindent\underline{\textbf{Step 4:} Applying $m \otimes 1$}

\noindent Finally, we apply $m$ to $\widetilde{B}_{j,q} \otimes A'_{u,s,j,q}$ and the identity to $A''_{u,s,j,q}$:
\begin{align}
\sum_{q,j,s,u}
\widetilde{B}_{j,q} \otimes A'_{u,s,j,q} \otimes A''_{u,s,j,q}
&\xmapsto{\;m \otimes 1\;}
\sum_{q=0}^{m} \sum_{j=1}^{j_q} \sum_{s=0}^{q} \sum_{u=1}^{u_s(j,q)}
\bigl(\widetilde{B}_{j,q} \times A'_{u,s,j,q}\bigr) \otimes A''_{u,s,j,q}. \notag 
\end{align}
%%%%%%%%%%%%%%%%%%%%%%%%%%%%%%%%%%%%%%%%%%%%%%%%%%%%%%%%%%%%%%%%%%%

Hence,
	$$M^*(Z( \Delta_1)\times\cdots\times Z( \Delta_r))=\sum_{q=0}^{m} \sum_{j=1}^{j_q} \sum_{s=0}^{q} \sum_{u=1}^{u_s(j,q)}\bigl(\widetilde{B}_{j,q} \times A'_{u,s,j,q}\bigr) \otimes A''_{u,s,j,q}. $$
	\noindent By Theorem \ref{main1}, $\mu^*(\pi\rtimes\sigma)=M^*(\pi)\rtimes \mu^*(\sigma),$  it follows that
	
	\begin{multline*}
		\mu^*(Z( \Delta_1)\times\cdots\times Z( \Delta_r)\rtimes \sigma)=
		\sum_{q=0}^{m}\sum_{j=1}^{j_q}\sum_{s=0}^{q}\sum_{u=1}^{u_s(j,q)} \Bigl( Z([-b_1, -{{p_1^{(j,q)}}}]_{(\tilde{\rho}_1)}) \times \cdots\times Z( [-b_r, -{{p_r^{(j,q)}}}]_{(\tilde{\rho}_r)}) \Bigr)\\
		\hspace{2.5cm}
		\times\Bigl( Z( [{a_1},{p_1^{(u_s(j,q))}-1}]_{(\rho_1)}) \times\cdots\times Z( [{a_r},{p_r^{(u_s(j,q))}-1}]_{(\rho_r)}) \Bigr)\Biggr]\\ 
		\otimes \Biggl[\Bigl( Z( [{p_1^{(u_s(j,q))}}, {p_1^{(j,q)}-1}]_{(\rho_1)}) \times\cdots\times Z( [{p_r^{(u_s(j,q))}},{p_r^{(j,q)}-1}]_{(\rho_r)}) \Bigr) \rtimes\sigma \Biggr].
	\end{multline*}
    
Since $\mu^*(Z(\Delta_1)\times\cdots\times Z(\Delta_r)\rtimes \sigma)$ is the sum of Jacquet modules with respect to all maximal Levi subgroups of $\Sp_n$, we fix $\M=\GL_l\times \Sp_{n-l}$ and retain only those summands of type $\sigma_1\otimes \sigma_2$ with $\sigma_1\in \mathcal{M}(\GL_l)$ and $\sigma_2\in \mathcal{M}(\Sp_{n-l})$. Thus
\begin{multline}\label{***}
r_{\M}^{\Sp_{n}}(Z(\Delta_1)\times\cdots\times Z(\Delta_r)\rtimes \sigma)=\\
\sum_{j=1}^{j_{l'}}\sum_{s=0}^{l'}\sum_{u=1}^{u_s(j,l')}
\Biggl[\Bigl( Z([-b_1,-p_1^{(j,l')}]_{(\tilde{\rho}_1)}) \times\cdots\times Z([-b_r,-p_r^{(j,l')}]_{(\tilde{\rho}_r)}) \Bigr)
\\
\hspace{1.5cm}\times\Bigl( Z([a_1,p_1^{(u_s(j,l'))}-1]_{(\rho_1)}) \times\cdots\times Z([a_r,p_r^{(u_s(j,l'))}-1]_{(\rho_r)}) \Bigr)\Biggr]
\\
\otimes \Biggl[\Bigl( Z([p_1^{(u_s(j,l'))},p_1^{(j,l')}-1]_{(\rho_1)}) \times\cdots\times Z([p_r^{(u_s(j,l'))},p_r^{(j,l')}-1]_{(\rho_r)}) \Bigr)\rtimes \sigma \Biggr].
\end{multline}
with $$
\sum_{i=1}^{r} m_i\bigl(b_i-p_i^{(j,l')}+1\bigr)=m-l',\qquad
\sum_{i=1}^{r} m_i\bigl(p_i^{(u_s(j,l'))}-a_i\bigr)=s,
$$

$$
\sum_{i=1}^{r} m_i\bigl(p_i^{(j,l')}-p_i^{(u_s(j,l'))}\bigr)=l'-s,\qquad
(m-l')+s=l.
$$

\begin{rem}\label{rem2} 
        We can observe that $r_{\M}^{\Sp_{n}}(Z( \Delta_1)\times\cdots\times Z( \Delta_r)\rtimes \sigma)$ has components of type $\sigma_i\otimes\pi_j$ (say), where $\sigma_i\in\mathcal{M}(\GL_l)$ and $\pi_j\in\mathcal{M}(\Sp_{n-l})$. If we take $Z( \Delta_1)\times\cdots\times Z( \Delta_r)\rtimes \sigma\in\mathcal{M}_{\Irr}(\Sp_{n})$, then by Lemma \ref{remark} and \ref{rem}, we conclude that $\sigma_i$ and $\pi_j$ are distinct irreducible representations of $\GL_l$ and $\Sp_{n-l}$, respectively.
	
\end{rem}

         Now, we prove that if the Jacquet functor $r_{\M}^{\Sp_n}$ is restricted to $\mathcal{M}_{\Irr}(\Sp_{n})$, then $$ \Dim_{\mathbb{C}}\Hom_{\M}(r_{\M}^{\Sp_{n}}(\pi), \rho) \leq 1,$$ for $\pi\in\mathcal{M}_{\Irr}(\Sp_{n})$ and any irreducible representation $\rho\in\mathcal{M}(\M)$.

	\subsection{Proof of Theorem \ref{main}.}
	\begin{proof}
		Let $m=(n_1,\ldots,n_r)$ be the partition of $m$ and $\Delta_i=[{a_i},{b_i}]_{(\rho_i)}$ be the non-linked  segments, where $\rho_i$ is an irreducible non self dual supercuspidal representation of $\GL_{m_i}$, for $i=1,2,\ldots,r.$  Let $b_i-a_i+1=k_i$ and $k_im_i=n_i$. This implies $Z( \Delta_1)\times\cdots\times Z( \Delta_r) $ is a representation of $\GL_m$.~Let $\pi=Z( \Delta_1)\times\cdots\times Z( \Delta_r)\rtimes \sigma\in\mathcal{M}_{\Irr}(\Sp_{n})$, where  $Z( \Delta_1)\times\cdots\times Z( \Delta_r)\in\mathcal{M}_{\Irr}(\GL_{m})$ and $\sigma\in\mathcal{C}_{\Irr}(\Sp_{n-m})$. Consider an irreducible representation $\rho=\tau_1\otimes\tau_2\in \mathcal{M}(\M)$. From (\ref{***}), it follows that
		\begin{multline*}
			r_{\M}^{\Sp_{n}}(\pi)=
			\sum_{j=1}^{j_{l'}}\sum_{s=0}^{l'}\sum_{u=1}^{u_s(j,l')}  \Biggl[\Bigl( Z( [-b_1, -{{p_1^{(j,l')}}}]_{(\tilde{\rho}_1)})\times\cdots\times Z( [-b_r, -{{p_r^{(j,l')}}}]_{(\tilde{\rho}_r)}) \Bigr)
			\\
			\hspace{1.5cm}	\times\Bigl( Z( [{a_1},{p_1^{(u_s(j,l'))}-1}]_{(\rho_1)}) \times\cdots\times Z( [{a_r},{p_r^{(u_s(j,l'))}-1}]_{(\rho_r)}) \Bigr)\Biggr]
			\\
			\otimes \Biggl[\Bigl( Z( [{p_1^{(u_s(j,l'))}}, {p_1^{(j,l')}-1}]_{(\rho_1)}) \times\cdots\times Z( [{p_r^{(u_s(j,l'))}},{p_r^{(j,l')}-1}]_{(\rho_r)}) \Bigr)\rtimes \sigma \Biggr].\hspace{0.5cm}
		\end{multline*} 
		 Observe that  $r_{\M}^{\Sp_{n}}(\pi)$ has subquotients of type $\sigma_i\otimes\pi_j$ (say). It is clear that $r_{\M}^{\Sp_{n}}(\pi)$ contain components of type $\widetilde{\sigma_i}\otimes\pi_j$ as well. Furthermore, if $Z(\Delta_1')\times\cdots\times Z(\Delta_t')\otimes \pi_j$ appears as a subquotient in $r_{\M}^{\Sp_{n}}(\pi)$, then none of its permutations are subquotients in  $r_{\M}^{\Sp_{n}}(\pi)$. Since, $\rho_i\ncong\widetilde{\rho_i}$, by  Remarks \ref{rem} and \ref{rem2}, all the components in $r_{\M}^{\Sp_{n}}(\pi)$ are distinct irreducible representations of $\GL_l\times\Sp_{n-l}$. By Theorem \ref{isom}, the representation $\tau_1\otimes\tau_2$ of $\M$ is either isomorphic to a component or does not appear in $r_{\M}^{\Sp_{n}}(\pi)$ at all. Hence, it follows that $\Dim_{\mathbb{C}}\Hom_{\M}(r_{\M}^{\Sp_{n}}(\pi), \rho)  
		\leq 1.$
		
	\end{proof}

  As observed above, the computation of Jacquet modules is quite complicated. To illustrate this complexity, we provide the following example.
	
	% \begin{example}
	% 	For $\Sp_2:$
	% 	Let $\rho$ be an irreducible representation of $\GL_1$ and $\sigma$ be an irreducible supercuspidal representation of $\Sp_1$. Then $\mu^*(\rho\rtimes \sigma)$ is computed as follows:
	% 	$$\rho$$
	% 	$$\bigg\downarrow s\circ m^*$$
	% 	$$\rho\otimes 1+1\otimes\rho$$
	% 	$$\bigg\downarrow \widetilde{.}\otimes m^*$$
	% 	$$\widetilde{\rho}\otimes 1\otimes 1+\widetilde{1}\otimes\rho\otimes 1+\widetilde{1}\otimes 1\otimes\rho$$
	% 	$$\bigg\downarrow m\otimes 1$$
	% 	$$\widetilde{\rho}\otimes 1+\rho\otimes 1+1\otimes\rho$$
	% 	Since $\mu^*(\sigma)=1\otimes\sigma,$
	% 	$\mu^*(\rho\rtimes\sigma)=\widetilde{\rho}\otimes \sigma+\rho\otimes \sigma+\rho\rtimes\sigma$.
	% 	\\
	% 	This yields:
	% 	\begin{itemize}
	% 		\item 
			
	% 		$r_{(1)}(\rho\rtimes\sigma)=\widetilde{\rho}\otimes \sigma+\rho\otimes \sigma.$
	% 		\item $r_{(2)}(\rho\rtimes\sigma)=0.$
	% 	\end{itemize}
	% 	It is clear that for $l=1~\text{and}~2,$ $$\Dim_{\mathbb{C}}\Hom_{\M}(r_{(l)}(\rho\rtimes\sigma), \tau)\leq1,$$  provided $\rho$ is non self dual.
	% \end{example}
	\begin{example}
	     Let $\rho_1,\rho_2\in \Irr(\GL_1)$ and $\sigma\in\mathcal{C}_{\Irr}(\Sp_1)$ be such that $\rho_1\times\rho_2\rtimes \sigma\in\Irr(\Sp_3).$ Using $\mu^*(\pi\rtimes\sigma)=M^*(\pi)\rtimes\mu^*(\sigma),$ we obtain
\begin{align*}
\mu^*(\rho_1\times\rho_2\rtimes\sigma)
&= \rho_1\times\rho_2\rtimes\sigma
+ \rho_1\otimes(\rho_2\rtimes\sigma)
+ \rho_2\otimes(\rho_1\rtimes\sigma) \\
&\quad + (\rho_1\times\rho_2)\otimes\sigma
+ \widetilde{\rho_1}\otimes(\rho_2\rtimes\sigma)
+ \widetilde{\rho_2}\otimes(\rho_1\rtimes\sigma) \\
&\quad + (\widetilde{\rho_1}\times\rho_2)\otimes\sigma
+ (\widetilde{\rho_2}\times\rho_1)\otimes\sigma
+ (\widetilde{\rho_1}\times\widetilde{\rho_2})\otimes\sigma.
\end{align*}
Hence,
\begin{align*}
r_{(1)} &= \rho_1\otimes(\rho_2\rtimes\sigma)
+ \rho_2\otimes(\rho_1\rtimes\sigma)
+ \widetilde{\rho_1}\otimes(\rho_2\rtimes\sigma)
+ \widetilde{\rho_2}\otimes(\rho_1\rtimes\sigma),\\
r_{(2)} &= (\rho_1\times\rho_2)\otimes\sigma
+ (\widetilde{\rho_1}\times\rho_2)\otimes\sigma
+ (\widetilde{\rho_2}\times\rho_1)\otimes\sigma
+ (\widetilde{\rho_1}\times\widetilde{\rho_2})\otimes\sigma,\\
r_{(3)} &= 0.
\end{align*}
Moreover, for $1\le l\le 3$ and any irreducible representation $\tau$ of maximal Levi $\M\simeq \GL_l\times\Sp_{(3-l)}$,
\[
\dim_{\mathbb{C}}\Hom_{\M}\!\big(r_{(l)}(\rho_1\times\rho_2\rtimes\sigma),\tau\big)\le 1,
\]
provided $\rho_1\not\simeq\rho_2$ and both are non-self-dual.

\end{example}

	\subsection{Acknowledgement:} The first author gratefully acknowledges the Ministry of Human Resource Development (MHRD), Govt. of India, for providing the necessary funding and fellowship to pursue this research work. 
	\section{Declaration}
	\textbf{Declaration for Research Article:}
	We, the authors of the research article ``On the Jacquet functor of Symplectic groups", declare that this article represents original work. It has not been published previously and is under consideration for publication elsewhere. We take full responsibility for the content of this article and affirm that it does not contain any defamatory or unlawful material.

	\bibliography{spjacq}
	\bibliographystyle{alpha}
\end{document}